\documentclass[12pt]{article}
\usepackage{amssymb,amsfonts}
\usepackage{color}
\usepackage{latexsym}
\usepackage{enumerate}
\usepackage{graphics}
 \usepackage{tikz}
\usepackage{graphicx}
\usepackage{mathdots}
\usepackage{bm}
\usepackage{amsmath}
\usepackage{mathrsfs}
\usepackage{pdfpages}

\newcommand\red[1] {{\color{red} #1}}
\newcommand\green[1] {{\color{green} #1}}
\newcommand\blue[1] {\textbf{ \color{blue} #1}}

\newtheorem{lem}{Lemma}
\newtheorem{theo}{Theorem}
\newtheorem{pro}{Proposition}
\newtheorem{prob}{Problem}
\newtheorem{cla}{Claim}
\newtheorem{prop}{Property}

\newcommand{\proof}{{\noindent {\em Proof}.\quad}\setcounter{countclaim}{0}
\setcounter{countcase}{0}}
\newcommand{\proofend}{{\hfill$\Box$}}

\newcounter{countfig}

\newcounter{countclaim}

\newcounter{countcase}

\newcommand{\beeq}{\begin{equation}}
\newcommand{\eneq}{\end{equation}}

\newcommand{\beeqn}{\begin{eqnarray*}}
\newcommand{\eneqn}{\end{eqnarray*}}

\def \usecolour
{\newcommand {\rered}{\red}
\newcommand {\reblue} {\blue}
\newcommand {\regreen} {\green}
}

\usecolour   
\def \R {{\cal R}}
\def \C {{\cal C}}
\def \al {\alpha}

\newcommand\equ[2]
{
\begin{equation}\label{#1}
	#2
\end{equation}
}

\newcommand\eqn[2]
{
	\begin{eqnarray} \label{#1}
		#2
	\end{eqnarray}
}

\usepackage{pdfpages}
\usepackage{tikz}
\usetikzlibrary{snakes}
\usetikzlibrary{backgrounds}
\usetikzlibrary{arrows}
\usetikzlibrary{shapes,shapes.geometric,shapes.misc}
\pgfdeclarelayer{edgelayer}
\pgfdeclarelayer{nodelayer}
\pgfsetlayers{background,edgelayer,nodelayer,main}
\tikzstyle{none}=[inner sep=0mm]

\tikzstyle{bluenode}=[fill=blue, draw=black, shape=circle, minimum size=0.15cm, inner sep=0pt]
\tikzstyle{blacknode}=[fill=black, draw=black, shape=circle, minimum size=0.15cm, inner sep=0pt]
\tikzstyle{pinknode}=[fill={rgb,255: red,255; green,191; blue,191}, draw=black,
shape=circle, minimum size=0.15cm, inner sep=0pt]
\tikzstyle{rednode}=[fill={rgb,255: red,244; green,0; blue,0}, draw=black,
shape=circle, minimum size=0.15cm, inner sep=0.1pt]
\tikzstyle{whitenode}=[fill={rgb,255: red,245; green,245; blue,245},
draw=black, shape=circle, minimum size=0.1cm, inner sep=0pt]
\tikzstyle{springnode}=[fill={rgb,255: red,44; green,218; blue,154},
draw=black, shape=circle, minimum size=0.15cm, inner sep=0pt]
\tikzstyle{star}=[fill={rgb,255: red,233; green,0; blue,81}, draw=red,
shape=star, minimum size=0.1cm, inner sep=0pt, tikzit fill={rgb,255: red,255;
green,191; blue,191}, tikzit draw={rgb,255: red,174; green,91; blue,32}]
\tikzstyle{orangenode}=[fill={rgb,255: red,255; green,128; blue,0}, draw=black,
shape=circle, minimum size=0.15cm, inner sep=0pt]
\tikzstyle{whilenode2}=[fill=white, draw=white, shape=circle]
\tikzstyle{testnode}=[fill=none, draw=black, shape=circle, minimum size=2.8cm,
inner sep=0.3pt]
\tikzstyle{testnode2}=[fill=none, draw=black, shape=circle, minimum size=0.2cm,
inner sep=0.1pt]
\tikzstyle{whitenode2}=[fill={rgb,255: red,242; green,242; blue,242},
draw=black, shape=circle, minimum size=0.12cm, inner sep=2pt]
\tikzstyle{green_dash}=[dash pattern=on 0.2mm off 0.2mm, draw=green, dashed]
\tikzstyle{dot}=[fill=black, draw=black, shape=circle, minimum size=0.04cm, inner sep=0pt]
\tikzstyle{blackedge}=[-, draw=black, fill=none, line width=0.22mm]
\tikzstyle{blackedge_thick}=[-, draw=black, line width=0.45mm, fill=none]
\tikzstyle{balck_dash}=[-, dash pattern=on 0.2mm off 0.2mm]
\tikzstyle{greenedge}=[-, draw={rgb,255: red,9; green,122; blue,43}]
\tikzstyle{pureedge}=[-, draw={rgb,255: red,230; green,0; blue,230}]
\tikzstyle{blue_directedge}=[draw=blue, ->]
\tikzstyle{rededge}=[-, draw=red]
\tikzstyle{rededge_thick}=[-, line width=0.45mm, draw=red]
\tikzstyle{red_dash}=[-, draw=red, dash pattern=on 0.2mm off 0.2mm]
\tikzstyle{dashpure1}=[-, dashed, draw={rgb,255: red,128; green,0; blue,128}]
\tikzstyle{blueedge}=[-, draw={rgb,255: red,0; green,0; blue,224}, line
width=0.22mm]
\tikzstyle{blue_dash}=[-, draw={rgb,255: red,13; green,20; blue,109}, line
width=0.1mm, dash pattern=on 0.2mm off 0.2mm]
\tikzstyle{dashedge}=[-, dash pattern=on 0.2mm off 0.2mm, draw={rgb,255:
red,143; green,0; blue,50}]
\tikzstyle{blackedge2}=[-, draw=black]
\tikzstyle{blackarrow}=[draw=black, ->]
\tikzstyle{orangeedge}=[-, draw={rgb,255: red,255; green,128; blue,0}]
\tikzstyle{whiteedge}=[-, draw=white, fill=none]
\tikzstyle{new edge style 0}=[-, draw=blue, line width=0.45mm]
\tikzstyle{new edge style 1}=[-, fill=white]
\tikzstyle{shadow_Lightgrey}=[-, fill={rgb,255: red,224; green,224; blue,224},
draw=black]
\tikzstyle{shadow_Lightgreen}=[-, draw=black, fill={rgb,255: red,204;
green,250; blue,223}]
\tikzstyle{shadow}=[-, fill={rgb,255: red,130; green,202; blue,255}, draw=red]
\tikzstyle{shadow_blue}=[-, fill={rgb,255: red,130; green,202; blue,255},
draw=black]
\tikzstyle{shadowblack}=[-, fill={rgb,255: red,240; green,248; blue,255},
draw=black]
\tikzstyle{shadow_pure}=[-, fill={rgb,255: red,224; green,182; blue,250},
draw=none]
\tikzstyle{shadow_}=[-, fill={rgb,255: red,145; green,145; blue,145}, draw=none]
\tikzstyle{shadow_silver}=[-, draw=black, fill={rgb,255: red,186; green,186;
blue,186}]
\tikzstyle{pickshing}=[-, draw={rgb,255: red,128; green,0; blue,128},
fill={rgb,255: red,255; green,191; blue,191}]

\begin{document}

\title {\textbf{The density of maximal IC-plane graphs and maximal NIC-plane graphs\thanks {The work was supported by the National
Natural Science Foundation of China (No.\ 12371346, 12271157, 12371340), Hunan Provincial Natural Science Foundation of
China (No.\ 2023JJ30178)
and the Scientific Research Fund of Hunan Provincial Education Department (No.\ 22A0637). Panna Geh\'er was supported by ERC Advanced Grant `GeoScape' No.\ 882971 and by the Thematic Excellence Program TKP2021-NKTA-62 of the National Research, Development and Innovation Office.}}}

\author
{Zongpeng Ding\thanks{Email: dzppxl@163.com}\\
\small School of Mathematics and Statistics\\
\small Hunan First Normal University,
Changsha 410205, China\\
\\
Yuanqiu Huang\thanks{Corresponding author. Email: hyqq@hunnu.edu.cn}\\
\small  School of Mathematics and Statistics\\
\small Hunan Normal University, Changsha 410081, China\\
\\
Fengming Dong\thanks{Email: donggraph@163.com
and fengming.dong@nie.edu.sg}\\
\small National Institute of Education,  Nanyang Technological
University,  Singapore\\
\\
Shengxiang Lv\thanks{ Email: lvsxx23@126.com}\\
\small  School of Mathematics and Statistics\\
\small Hunan University of Finance and Economics, Changsha 410205, China\\
\\
Panna Geh\'{e}r\thanks{ Email: geher.panna@ttk.elte.hu}\\
\small  E\"otv\"os University, Hungary\\
\small Alfr\'ed R\'enyi Institute of Mathematics, Hungary\\
\\
}

\date{}

\maketitle

\newpage

\begin{abstract}
In this paper, we show that
any maximal IC-plane graph of order $n$ has at least $\left\lceil\frac{7}{3}n-\frac{14}{3}\right\rceil$ edges,
and any maximal NIC-plane graph of order $n$ has at least $\left\lceil\frac{11}{5}n-\frac{18}{5}\right\rceil$ edges.
Moreover, we show that both results are tight for infinitely many integers $n$.
\end{abstract}




\section{Introduction}
Geometric graph theory investigates the combinatorial properties of  graphs drawn in the Euclidean plane with possibly intersecting edges.
A \textit{drawing} of a finite, undirected graph $G=(V, \, E)$ is
a mapping $D$ that assigns each vertex $v \in V(G)$ a distinct point $D(v)$ in the plane, and assigns each edge $uv \in E(G)$
a continuous arc connecting the points $D(u)$ and $D(v)$, not passing through the image of any other vertex. For terminology and notation not defined here, we refer to \cite{BA}.
A graph $G$ is \emph{maximal} in a given graph class if no edge can be added to $G$ without violating any condition of that class. Determining the density of maximal graphs in a given graph class is a central problem in geometric graph theory.
Recall that  due to the famous Euler's formula, a planar graph with $n\geq 3$ vertices has at most $3n-6$ edges. Moreover, any maximal planar graph with $n\geq 3$ vertices, is a triangulation.

In recent research topics, several relaxations of planarity have been considered. In this context, we allow some crossings, but we forbid specific edge crossing patterns.
A graph is called \textit{$1$-planar} if it can be drawn in the plane such that
each edge is crossed at most once. The notion of 1-planarity was introduced by Ringel~\cite{RG} in 1965, and since then many properties of 1-planar graphs have been investigated (see, e.g.~\cite{AE,BM,BJ,CH,CJ,KV,SY}). The density of this class is also well-studied. It is known that a 1-planar graph has at most $4n-8$ edges (see~\cite{BR,PJ}) and this bound is tight. However, there is a fundamental difference in the density of 1-planar and planar graphs. Unlike planar graphs, maximal $1$-planar graphs can have different densities. Thus, it is interesting to study the minimum number of edges a maximal 1-planar graph can have. As already mentioned, on $n$ vertices there exists a maximal $1$-planar graph with $4n-8$ edges (see~\cite{BR}), but also there exists one with as few as $\frac{45}{17}n-\frac{84}{17}$ edges (see~\cite{BF}). In general, Bar\'at and T\'oth~\cite{BJ} proved that any maximal 1-planar graph of order $n$ has at least $\frac{20}{9} n - \frac{10}{3}$ edges. Recently, this result was improved by Huang et~al.~\cite{H}. They showed that any maximal $1$-planar graph of order $n \geq 5$ contains at least
$\left \lceil \frac {7}3n \right \rceil-3$ edges.

Subclasses of $1$-planar graphs also have rich literature. A graph is called \textit{IC-planar} (where IC stands for independent crossings) if it admits a drawing such that each edge is crossed at most once and any two pairs of crossing edges share no common end vertex. This graph class was introduced by Albertson~\cite{AM} (see also~\cite{BFJ,KD}). A graph is called \textit{NIC-planar} (where NIC stands for near-independent crossings) if it admits a drawing such that each edge is crossed at most once and any two pairs of crossing edges share at most one vertex. This graph class was introduced by Zhang~\cite{Z}
(see also~\cite{B}).

The density of these graph classes is also well-studied.
For IC-planar graphs, Zhang and Liu~\cite{ZX} showed that maximal IC-planar graphs of order $n$ have at most $\frac{13}{4}n-6$ edges, and this bound is tight.
Ouyang et~al.~\cite{O} showed that an IC-planar drawing of any maximal IC-planar graph of order $n\ge 5$ is a triangulation. Using this result, it is not difficult to prove that maximal IC-planar graphs
of order $n\ge 5$ have at least $3n-5$ edges, and this bound is tight, see~\cite{BC}.
For NIC-planar graphs, Zhang~\cite{Z} (and later also Bachmaier et~al.~\cite{B}) proved that a NIC-planar graph of order $n$ has at most $\frac{18}{5}n-\frac{36}{5}$ edges.
Czap and \v{S}ugerek~\cite{CP} found a NIC-planar graph with 27 vertices and 90 edges, which shows that the bound is tight for $n=27$.
For the lower bound, Bachmaier et~al.~\cite{B} showed that every maximal NIC-planar graph
of order $n\geq5$ has at least $\frac{16}{5}n-\frac{32}{5}$ edges. This result is tight for infinitely many values of~$n$.
A special case was also studied.
A vertex is called a \emph{dominating vertex} if it is adjacent to all other vertices. Ouyang et~al.~\cite{OZ} proved that a NIC-planar graph of order~$n$ with a dominating vertex has at most $\frac{7}{2}n-\frac{15}{2}$ edges, and this bound is tight.

The known results on the densities of maximal 1-planar graphs and its subclasses, IC-planar and NIC-planar graphs are summarized in Table~\ref{biaob1}.

\begin{table}[htbp]
\centering
\begin{tabular}{|c|c|c|c|}
  \hline
 & 1-planar & NIC-planar & IC-planar
 \\
\hline
\rule{0pt}{15pt}%
\rule[-1.5ex]{0pt}{0pt}%
Upper bound & $4n-8$ (see \cite{BR, PJ}) & $\frac{18}{5}n-\frac{36}{5}$ (see \cite{B, Z}) & $\frac{13}{4}n-6$ (see \cite{ZX})
\\  \hline
\rule{0pt}{15pt}%
\rule[-1.5ex]{0pt}{0pt}%
Lower bound &$\left \lceil \frac {7}3n \right \rceil-3$  (see \cite{H})  & $\frac{16}{5}n-\frac{32}{5}$ (see \cite{B}) & $3n-5$ (see \cite{BC})
\\  \hline
\end{tabular}
\caption{The density of maximal $1$-planar, NIC-planar and IC-planar graphs.} \label{biaob1}
\end{table}

In this paper, we examine graphs with given drawings. We say that a drawing of a graph is \emph{saturated} in a graph class if it satisfies all conditions of the class, but no more edges can be added to the drawing without violating one of them.
Since any planar graph with a planar drawing (that is called a  plane graph) can be extended to a triangulation, maximal plane graphs are exactly the maximal planar graphs. In the spirit of the notion of plane graphs, a graph together with a
$1$-planar drawing is called
a $1$-plane graph. IC-plane and NIC-plane graphs can be defined analogously. Furthermore, we say that a 1-plane (resp., IC-plane, NIC-plane) graph with a saturated drawing is a maximal 1-plane (resp., IC-plane, NIC-plane) graph.



\begin{figure}[h]
\setlength{\unitlength}{1.5mm}
     \begin{center}
  \begin{tikzpicture}[scale=0.8]

\coordinate(u1) at  (1, 1);
\coordinate(v1) at  (1, 4);
\coordinate(u2) at  (4, 4);
\coordinate(v2) at  (4, 1);

 \coordinate(c1) at  (1.5, 2.5);

\draw[thick](u1)--(u2)(v1)--(v2)(u1)--(v1)(v1)--(u2)(u2)--(v2)(v2)--(u1);
\draw[thick](u1)--(c1)(v1)--(c1);

 \filldraw [black](u1) circle (3pt)(u2) circle (3pt)(v1) circle (3pt)(v2) circle(3pt);
\filldraw [black](c1) circle (3pt);
  \node at (2.5, 0.4) {(a)};

 \coordinate(u21) at  (9, 1);
\coordinate(v21) at  (9, 4);
\coordinate(u22) at  (12, 4);
\coordinate(v22) at  (12, 1);

 \coordinate(c21) at  (8.5, 2.5);

\draw[thick](u21)--(u22)(v21)--(v22)(u21)--(v21)(v21)--(u22)(u22)--(v22)(v22)--(u21);
\draw[thick](u21)--(c21)(v21)--(c21);

 \filldraw [black](u21) circle (3pt)(u22) circle (3pt)(v21) circle (3pt)(v22) circle(3pt);
\filldraw [black](c21) circle (3pt);
  \node at (10.5, 0.4) {(b)};

     \end{tikzpicture}
  \caption{{\small Two different IC-planar drawings of an IC-planar graph $G$. The drawing on the left is saturated, while drawing on the right is not. Hence, $G$ is not a maximal IC-planar graph.}}\label{1}
  \end{center}

  \end{figure}
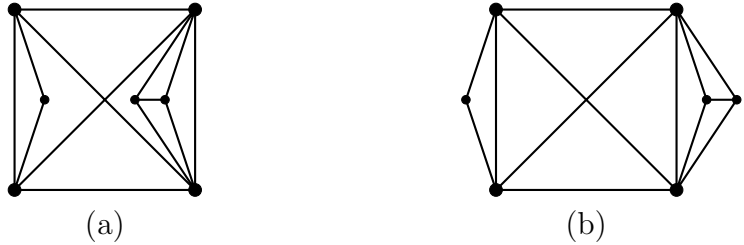

The above mentioned upper bounds for the edge number of 1-planar, IC-planar and NIC-planar graphs also give tight upper bounds for 1-plane, IC-plane and NIC-plane graphs. However, a maximal 1-plane graph might have less edges than a maximal 1-planar graph has with the same number of vertices. This is because for a maximal 1-plane graph, the underlying 1-planar graph may not be maximal. This observation also holds for IC-plane and NIC-plane graphs, an example is shown in Figure~\ref{1}.

Thus, the following problem arises.
\begin{prob}
    Determine the minimum numbers of edges of maximal 1-plane, IC-plane and NIC-plane graphs, respectively.
\end{prob}



As we have already mentioned,
Huang et~al.~\cite{H} showed that every maximal
$1$-planar graph of order $n \geq 5$ has at least
$\left \lceil \frac {7}3n \right \rceil-3$ edges. In fact, they proved a stronger statement that their result also holds for maximal
$1$-plane graphs. This bound is tight for all $n\ge 5$. Recently, Xu~\cite{X} showed that any maximal IC-plane graph with $n \geq 4$ vertices has at least $2n-2$ edges. However, the author mentioned that this result seems to be far from optimal and called the problem of determining the exact value an intriguing question. In this paper, we solve this problem and we also derive a tight bound on the minimum density for NIC-plane graphs. We have the following results.

\begin{theo}\label{theo2}
Any maximal IC-plane graph of order $n$
has at least $\left\lceil\frac{7}{3}n-\frac{14}{3}\right\rceil$ edges.
Moreover, this bound is tight for
all integers $n\ge 8$
with $n\equiv 2\pmod{6}$.
\end{theo}

\begin{theo}\label{theo2N}
Any maximal NIC-plane graph of order $n$
has at least $\left\lceil\frac{11}{5}n-\frac{18}{5}\right\rceil$ edges.
Moreover, this bound is tight for
all integers $n\ge 8$
with $n\equiv 3\pmod{5}$.
\end{theo}

\noindent \textbf{Paper outline.}
The aim of Section~\ref{preliminaries} is to introduce some basic properties of maximal IC-plane and maximal NIC-plane graphs. In Section \ref{tightness}, we construct maximal IC-plane and maximal NIC-plane graphs that demonstrate the tightness of our main theorems. Finally, in Section~\ref{proof}, we prove our main results.

\section{Preliminaries} \label{preliminaries}
Throughout this section, we always assume that $G$ is a maximal IC-plane graph or a maximal NIC-plane graph of order $n$ with $n\ge 4$.
If two edges $ab$ and $cd$ in $E(G)$
cross each other at a point $\alpha$,
then they are called \textit{crossing
	edges},
and $\alpha a$, $\alpha b$,
 $\alpha c$ and $\alpha d$ are
 called \textit{half-edges}.
 A \textit{non-crossing} edge in $G$ does not
 cross any other edge in $G$.
The edges in $G$ divide the plane into
\textit{faces}, each of which is bounded by edges and half-edges.
\textit{The size} of a face refers to the number of edges and half-edges on the boundary of the face.
A face is called a \textit{$k$-face} if its size
is $k$.
A face of $G$ is called
a \textit{true face} if
it is bounded by non-crossing edges only,
and called a \textit{false} face otherwise.
For any face $F$ of $G$, let $\partial(F)$
denote its boundary.
Some basic properties on $G$ are introduced below.

\begin{prop}\label{propa}
Let $G$ be either a maximal IC-plane or a maximal NIC-plane graph and $F$ be a face in $G$.
\begin{enumerate}[(1)]
\item Any two vertices
on $\partial(F)$ are adjacent in $G$.

\item
$\partial(F)$ contains at least two and at most four vertices.
Moreover, if $\partial(F)$ has exactly two vertices, then $F$ is a false $3$-face;
and
if $\partial(F)$ has four vertices,
then $F$ is a true $4$-face.

\item If edges $ab$ and $cd$ in $G$ cross
each other,  then
$\{a, \, b, \,  c, \, d\}$ form a clique in $G$.
	
\item The minimum degree of $G$
is at least $2$.
\end{enumerate}
\end{prop}

These properties can be checked analogously to those in maximal $1$-plane graphs (see \cite[Lemma 1]{BF} and \cite[Lemmas 1, 2 and 3]{BJ}). Thus, we omit the proof.

The following lemma describes the possible structures of faces in maximal IC-plane or maximal NIC-plane graphs.

\begin{figure}[h]
	\setlength{\unitlength}{1.5mm}
	\begin{center}
		\begin{tikzpicture}[scale=0.8]
			
			\coordinate[label=left:$a$](a1) at  (0.5, 1);
			\coordinate[label=right:$b$](b1) at  (2, 1);
			\coordinate[label=left:$c$](c1) at  (1.25, 2.5);
			
			\draw[thick](a1)--(b1)--(c1)--(a1);
			\filldraw [black](a1) circle (3pt)(b1) circle (3pt)(c1) circle (3pt);
			\node at (1.2, 1.5) {$F$};
			\node at (1.25, 0.4) {(i)};

			\coordinate[label=left:$a$](a21) at  (3.5, 1);
			\coordinate[label=right:$b$](b21) at  (5, 1);
			\coordinate[](c21) at  (4, 2.7);
			\coordinate[](c22) at  (4.5, 2.7);
			
			\draw[thick](a21)--(b21)(a21)--(c22)(b21)--(c21);
			\filldraw [black](a21) circle (3pt)(b21) circle (3pt);
            \filldraw [black](c22) circle (3pt);
            \filldraw [black](c21) circle (3pt);
			\node at (4.2, 1.5) {$F$};
			\node at (4.2, 0.4) {(ii)};
			
			\coordinate[label=left:$a$](a31) at  (6.5, 1);
			\coordinate[label=right:$b$](b31) at  (8, 1);
			\coordinate[label=right:$c$](c31) at  (8, 2.7);
			\coordinate[label=left:$d$](d31) at  (6.5, 2.7);

			\draw[thick](a31)--(b31)--(c31)--(d31)--(a31);
			\draw[thick] (a31)..controls(5.6, 3.5)..(c31);
			\draw[thick] (b31)..controls(8.8, 3.5)..(d31);
			\node at (7.25, 3.2) {$\alpha$};
			
			\filldraw [black](a31) circle (3pt)(b31) circle (3pt)(c31) circle (3pt)(d31) circle (3pt);
			
			\node at (7.2, 1.7) {$F$};
			\node at (7.2, 0.4) {(iii)};
			
			
			\coordinate[label=left:$a$](a41) at  (9.5, 2);
			\coordinate[label=right:$b$](b41) at  (11, 2);
			\coordinate[label=right:$c$](c41) at  (11, 3);
			\coordinate[label=left:$d$](d41) at  (9.5, 3);
			\coordinate[label=left:$e$](e) at  (10.25, 1);

			\draw[thick](a41)--(c41)(d41)--(b41)(a41)--(e)(b41)--(e);
			\node at (10.25, 2.8) {$\alpha$};
			\draw[thick] (a41)..controls(8.6, 3.5)..(10.25, 3.5)..controls(11.9, 3.5)..(b41);
			
			\filldraw [black](a41) circle (3pt)(b41) circle (3pt)(c41) circle (3pt)(d41) circle (3pt)(e) circle (3pt);
			
			\node at (10.3, 1.7) {$F$};
			\node at (10.2, 0.4) {(iv)};
			
			\coordinate[label=left:$a$](a51) at  (12.6, 2.25);
			\coordinate[label=right:$b$](b51) at  (14, 2.25);
			\coordinate[label=right:$c$](c51) at  (14, 3.15);
			\coordinate[label=left:$d$](d51) at  (12.5, 3.15);
			\coordinate[label=left:$e$](e) at  (13.25, 1);
			\coordinate(e1) at  (12.5, 1.25);
			\coordinate(e2) at  (12.2, 1.85);
			\draw[thick](a51)--(c51)(d51)--(b51)(b51)--(e)(a51)--(e1)(e)--(e2);
			\node at (13.25, 2.95) {$\alpha$};
			\filldraw [black](a51) circle (3pt)(b51) circle (3pt)(c51) circle (3pt)(d51) circle (3pt)(e) circle (3pt)(e1) circle (3pt) (e2) circle (3pt);
			
			\node at (13.3, 1.85) {$F$};
			\node at (13.2, 0.5) {(v)};
			
			\coordinate[label=left:$a$](a61) at  (15.6,2.25);
			\coordinate[label=right:$b$](b61) at  (17, 2.25);
			\coordinate[label=right:$c$](c61) at  (17, 3.15);
			\coordinate[label=left:$d$](d61) at  (15.5,3.15);
			\coordinate[label=left:$e$](e) at  (16.25, 1);
			\coordinate(e1) at  (15.5, 1.15);
			\coordinate(e2) at  (15.2, 1.75);
			\coordinate(e3) at  (16.8, 1.15);
			\coordinate(e4) at  (17.2, 1.75);
			\draw[thick](a61)--(c61)(d61)--(b61)(a61)--(e1)(e)--(e2)(b61)--(e3)(e)--(e4);
			\node at (16.25, 2.95) {$\alpha$};
			
			\filldraw [black](a61) circle (3pt)(b61) circle (3pt)(c61) circle (3pt)(d61) circle (3pt)(e) circle (3pt)(e1) circle (3pt) (e2) circle (3pt)(e3) circle (3pt) (e4) circle (3pt);
			
			\node at (16.3, 1.85) {$F$};
			\node at (16.2, 0.5) {(vi)};
		\end{tikzpicture}
\caption{Possible types of faces in $G$.}
\label{figaN}
\end{center}
\end{figure}
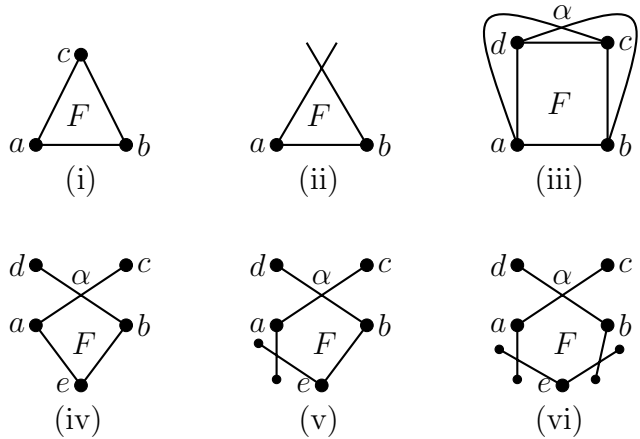

\begin{lem}\label{lemaN}
Assume that $G$ is either a maximal IC-plane graph or a maximal NIC-plane graph, and
$F$ is a $k$-face in $G$.
Then $3\leq k\leq 4$ if $G$ is
a maximal IC-plane graph, and
$3\leq k\leq 6$ otherwise.
Furthermore,
\vspace{-4 mm}
\begin{enumerate}[(1)]
		
	\item if $k=3$, then
	$F$ is either a true
	$3$-face or a false $3$-face,
	as shown in Figure \ref{figaN} (i) or (ii);

	\item  if $k=4$, then
	$F$ is either a true
	$4$-face or a false $4$-face,
	as shown in Figure \ref{figaN} (iii) or (iv); and

	\item  if $5\le k\le 6$,
then $F$ is a false $k$-face,
as shown in Figure \ref{figaN}
(v) or (vi).


\end{enumerate}
\end{lem}


\proof Clearly, since $F$ is a $k$-face, $k\ge 3$.
Suppose that $G$ is a maximal IC-plane graph  and $k\ge 5$. By Property~\ref{propa}~(2),
there are exactly three vertices
on $\partial(F)$.
On the other hand, as $G$ is IC-plane,
$\partial(F)$ contains at most
two half-edges, implying that $\partial(F)$
contains more than three vertices,
a contradiction.

Now assume that $G$ is a maximal NIC-plane graph and $k\ge 7$.
By Property~\ref{propa}~(2),
there are exactly three vertices
on $\partial(F)$.
As each edge in $G$ can be crossed
at most once,
$\partial(F)$ contains at most
six edges and half-edges.
Thus, $k\le 6$, a contradiction.
Hence, $3\le k\le 4$ if $G$ is maximal IC-plane, and $3\le k\le 6$
otherwise.

Finally, note that by Property~\ref{propa}~(2),
conclusions
(1), (2) and (3) follow directly.
\proofend

Let $\alpha$ be the crossing
of two edges $ac$ and $bd$ in $G$.
We say $\alpha$
and a face $F$
are \textit{incident},
denoted by $\alpha \sim F$  or $ F \sim \alpha $,
if either $a, \, b, \, c, \, d$ are all on $\partial(F)$,
or $\alpha$ is on $\partial(F)$.
For example, $\alpha \sim F$ happens in 
Figure~\ref{figaN} (iii), (iv), (v) and (vi), respectively.
Note that $\alpha$ lies on
$\partial(F)$ if and only if
$\partial(F)$ contains one vertex out of $\{a, \, c\}$ and one vertex out of $\{b, \, d\}$.

Below we list some properties of the incidence relation between
a crossing $\alpha$ and a face $F$
in $G$.

\begin{prop}\label{propb}
Let $G$ be either a
maximal IC-plane graph or a
maximal NIC-plane graph.
	\vspace{-4 mm}
\begin{enumerate}[(1)]
	\item Any $k$-face in $G$, where $4\le k\le 6$, is incident with exactly
	$k-3$ crossings.

\item Each crossing in $G$ can be incident with at most one true $4$-face,
and at most four false faces
of size at least $4$.
\end{enumerate}
\end{prop}

\proof
(1). First assume that $F$ is a true face.
Then, $F$ is a $4$-face. Let the vertices of $F$ be $a, \, b, \, c, \, d$, in this order, as shown in Figure~\ref{figaN} (iii).
By Property~\ref{propa}~(1),
$ac$ and $bd$ cross each other
at a point $\alpha$
outside of $F$. By definition,
$F\sim \alpha $.
Clearly, $F$ cannot be incident with any other crossings, thus, $F$ is incident with exactly one crossing.

Now, let $F$ be a false $k$-face
(see Figure~\ref{figaN} (iv), (v), and (vi)),
where $4\le k\le 6$.
Clearly, there are exactly $k-3$
crossings on $\partial(F)$,
and each of them is incident with $F$.
Assume that $F$ is incident with
another crossing $\alpha^{\prime}$
which is not on
$\partial(F)$. Let $\alpha^{\prime}$ be the crossing of two edges $a^{\prime}c^{\prime}$ and $b^{\prime}d^{\prime}$ in $G$. By definition,
 $a^{\prime}, \, b^{\prime}, \, c^{\prime}, \, d^{\prime}$ are all on $\partial(F)$,
namely, $F$ is a true 4-face, a contradiction.
Hence,
$F$ is incident with exactly
$k-3$ crossings.

(2). Let $\alpha$ be a crossing of the edges $ac$ and $bd$.
First, assume that $F$ is a true face. Then, $F$ must be the face bounded by the four non-crossing edges of $K_4$. This face is unique for any crossing, thus, the $F$ is the only true 4-face that is incident with $\alpha$.

Next, assume that $F$ is a false face. The vertices $\{a, \, b, \, c, \, d \}$ cannot all lie on the boundary of $F$, thus $\alpha$ must be on the boundary of $F$. Clearly, each crossing
can lie on the boundary of at most four faces. Thus, the conclusion holds.
\proofend

Assume that edges $ab$ and $cd$ cross each other at a crossing $\alpha$.
By Property~\ref{propa} (3),
$\{a, \, b, \,  c, \, d\}$
induces a $K_4$.
Clearly,  this is the only $K_4$
in which two edges cross each other
at $\alpha$.
For convenience, we say
this $K_4$ is \textit{spanned} by
$\alpha$, denoted by $K_4(\alpha)$.
$K_4(\alpha)$
divides the plane into five regions $R_{i}$ ($1\leq i\leq5$),
where each of
	$R_1, \, R_2, \, R_3$ and $R_4$
is bounded by one edge and two
half-edges,
as shown in
Figure~\ref{5} (a) or (b).
Clearly, $R_5$ is the only region  bounded by four edges.
Let ${\cal R}(\alpha)$ denote the
set $\{R_i: 1\le i\le 5\}$.
 \begin{figure}[h]
\setlength{\unitlength}{1.5mm}
     \begin{center}
  \begin{tikzpicture}[scale=0.8]

\coordinate(a1) at  (0.5, 1);
\coordinate(b1) at  (2.5, 1);
\coordinate(c1) at  (2.5, 3);
\coordinate(d1) at  (0.5, 3);
\draw[thick](a1)--(b1)--(c1)--(d1)--(a1)--(c1)(b1)--(d1);
 \filldraw [black](a1) circle (3pt)(b1) circle (3pt)(c1) circle (3pt)(d1) circle (3pt);
 \node at (1.5, 1.5) {\small{$R_1$}};   \node at (2.1, 2) {\small{$R_2$}};   \node at (1.5, 2.5) {\small{$R_3$}};   \node at (0.9, 2) {\small{$R_4$}};   \node at (0.2, 2) {\small{$R_5$}};


\coordinate(u1) at  (8.5, 1);
\coordinate(v1) at  (10.5, 1);
\coordinate(u2) at  (10.5, 3);
\coordinate(v2) at  (8.5, 3);
\draw[thick](u1)--(v1)--(u2)--(v2)--(v1)(u1)--(u2);
\draw[thick](v2)..controls(11.5,4)..(11.5, 2)..controls(11.5,0)..(u1);
 \filldraw [black](u1) circle (3pt)(v1) circle (3pt)(u2) circle (3pt)(v2) circle (3pt);

  \node at (9.5, 1.5) {\small{$R_1$}};
  \node at (10.1, 2) {\small{$R_2$}};
  \node at (9.5, 2.5) {\small{$R_3$}};
  \node at (8.9, 2) {\small{$R_4$}};
  \node at (11.0, 2) {\small{$R_5$}};
   \end{tikzpicture}

(a) \hspace{5.5 cm} (b)

  \caption{The region $R_{5}$ on the left is infinite, and on the right is finite.}\label{5}
  \end{center}
\end{figure}


Following \cite{BF}, a vertex of degree $2$ in a graph is called a \textit{hermit}. The following observations will play an important role in our main results.

\begin{lem}\label{lemc}
Assume $G$ is a maximal IC-plane graph. If $G$  has exactly $c$ crossings and $h$ hermits, then $c\leq \frac{n-h}{4}.$ 
\end{lem}
\proof
By Property~\ref{propa} (3), each crossing spans a $K_{4}$ in $G$, and any two $K_{4}$'s
spanned by two crossings
share no common end vertex by the assumption on $G$.
Since each hermit is of degree $2$,
it cannot be a vertex on any $K_4$
spanned by a crossing.
Thus,
$$
n=|V(G)|\ge h+ 4c,
$$
implying that
 $c\leq \frac{n-h}{4}$.
Thus, Lemma~\ref{lemc} holds.
\proofend

\begin{lem}\label{lemcN}
Assume that $G$ is a maximal NIC-plane graph. If $G$ has exactly $c$ crossings and $h$ hermits, then $c\leq \frac{|E(G)|-2h}{6}.$ 
\end{lem}
\proof
By Property~\ref{propa} (3), each crossing spans a $K_{4}$ in $G$, and any two $K_{4}$'s
spanned by two crossings
do not share any edge by the assumption on $G$.
Since each hermit is of degree $2$,
it cannot be a vertex on any $K_4$
spanned by a crossing.
Thus,
$$
|E(G)|\ge 2h+ 6c,
$$
implying that
 $c\leq \frac{|E(G)|-2h}{6}$.
Thus, Lemma~\ref{lemcN} holds.
\proofend

\begin{lem}\label{lemB}
For any (not necessarily simple) graph $H$ of order $n$,  if
there exists a drawing $D$ for $H$
such that
each face in $D$ is a $3$-face and $D$ has exactly $c$ crossings, then $|E(H)|=3n-6+c$.
\end{lem}

\proof
 Let $D^{P}$
 be the planarization of $D$ obtained by adding new vertices at each crossing point of $D$. Since every face of $D$ is a 3-face, $D^{P}$
is a maximal plane graph
and thus
$|E(D^{P})|=3|V(D^{P})|-6$.
Note that $|V(D^{P})|=|V(D)|+c$ and $|E(D^{P})|=|E(D)|+2c$.
Thus, $|E(H)|=|E(D)|=3n-6+c$.
\proofend

\section{Tightness of bounds of Theorems~\ref{theo2} and~\ref{theo2N}}
\label{tightness}




In this section, we show that
the lower bound in Theorems~\ref{theo2} is tight
for all $n\ge 8$ with
$n\equiv 2 \pmod{6}$,
and the lower bound in Theorems~\ref{theo2N} is tight
for all $n\ge 8$ with
$n\equiv 3 \pmod{5}$.

Let $H^*$ and $H_{1}$ be the IC-plane graphs shown in Figure \ref{4} (a) and (b), respectively.
For $k\geq2$,  $H_{k}$ is
obtained from $H_{k-1}$
by inserting $H^*$ into
a false 3-face $F$ of $H_{k-1}$
and identifying edge $uv$ in $H^*$
with the only non-crossing edge on $\partial(F)$.
As an example,  $H_2$ is shown in
Figure \ref{4} (c).


 \begin{figure}[h]
\setlength{\unitlength}{1.5mm}
     \begin{center}
  \begin{tikzpicture}[scale=0.8]

\coordinate(a1) at  (0.8, 2);
\coordinate(b1) at  (2.2, 2);
\coordinate(c1) at  (2.2, 3);
\coordinate(d1) at  (0.8, 3);
\coordinate[label=left:$u$](e1) at  (0.0, 1);
\coordinate[label=right:$v$](f1) at  (3, 1);

\draw[thick](e1)--(a1)--(b1)--(c1)--(d1)--(a1)(b1)--(f1)(e1)--(f1)(e1)--(b1);
\draw[thick](a1)..controls(0.1, 3.4)..(1.5, 3.7)..controls(2.9, 3.8)..(f1);
\draw[thick](a1)..controls(0.4, 3.5)..(c1);
\draw[thick](b1)..controls(2.7, 3.5)..(d1);

 \filldraw [black](a1) circle (3pt)(b1) circle (3pt)(c1) circle (3pt)(d1) circle (3pt)(e1) circle (3pt)(f1) circle (3pt);
 \node at (1.45, 0.2) {(a)~$H^*$};


\coordinate(a1) at  (5.5, 1);
\coordinate(b1) at  (8, 1);
\coordinate(c1) at  (8, 4);
\coordinate(d1) at  (5.5, 4);

\draw[thick](a1)--(b1)--(c1)--(d1)--(a1);
\draw[thick](a1)..controls(6, 3)..(c1);
\draw[thick](b1)..controls(7.5, 3)..(d1);

 \filldraw [black](a1) circle (3pt)(b1) circle (3pt)(c1) circle (3pt)(d1) circle (3pt);
  \node at (6.6, 0.2) {(b)~$H_1$};

\coordinate(a1) at  (11.55, 1.5);
\coordinate(b1) at  (12.45, 1.5);
\coordinate(c1) at  (12.45, 2);
\coordinate(d1) at  (11.55, 2);
\coordinate(e1) at  (10.5, 1);
\coordinate(f1) at  (13.5, 1);
\coordinate(a2) at  (10.5, 5);
\coordinate(b2) at  (13.5, 5);
\draw[thick](e1)--(a1)--(b1)--(c1)--(d1)--(a1)(b1)--(f1)(e1)--(f1)(e1)--(b1)(e1)--(a2)--(b2)--(f1);
\draw[thick](a1)..controls(10.7, 2.7)..(12, 2.8)..controls(12.8, 2.8)..(f1);
\draw[thick](a1)..controls(11.01, 2.7)..(c1);
\draw[thick](b1)..controls(12.99, 2.7)..(d1);
\draw[thick](e1)..controls(11.01, 3.37)..(b2);
\draw[thick](f1)..controls(13.12, 3.37)..(a2);
 \filldraw [black](a1) circle (3pt)(b1) circle (3pt)(c1) circle (3pt)(d1) circle (3pt)(e1) circle (3pt)(f1) circle (3pt)(a2) circle (3pt)(b2) circle (3pt);
  \node at (11.95, 0.2)  {(c)~$H_2$};

      \end{tikzpicture}
  \caption{$H^{*}$, $H_{1}$ and $H_{2}$.}\label{4}
  \end{center}
\end{figure}

Observe that for each $k\ge 1$,
	$H_{k}$ is a maximal IC-plane graph without hermits, and
$$
|V(H_k)|=4k,
\qquad
|E(H_{k})|=6+10(k-1)=10k-4.
$$
Also note that each $H_k$
has exactly $3+2(k-2)+3=2k+2$ false $3$-faces.
For each $k\ge 1$,
let $H'_k$ be the graph obtained
by inserting a hermit within each
false 3-face of $H_{k}$.

\begin{pro}\label{H'_k-lem}
For each $k\ge 1$, $H'_k$
is a maximal IC-plane graph
with $n$ vertices and
$\frac{7}{3}n-\frac{14}{3}$ edges, where $n=6k+2$.
\end{pro}
\proof Since $H_k$ has exactly $2k+2$
false $3$-faces, by the definition
of $H'_k$,
$$
n=|V(H_k)|+2k+2
=4k+2k+2=6k+2
$$
and
$$
|E(H'_k)|=|E(H_k)|+2(2k+2)=
10k-4+2(2k+2)=14k=
\frac{7}{3}n-\frac{14}{3}.
$$
\proofend


This shows the tightness of bound of Theorem \ref{theo2}.

Now, let $M^*$ and $M_{1}$ be the NIC-plane graphs shown in Figure \ref{fig3N} (a) and (b) respectively.
For any $k\geq2$,  let $M_{k}$ be
obtained from $M_{k-1}$
by inserting $M^*$ into
a false 3-face $F$ of $M_{k-1}$
and identifying edge $uv$ in $M^*$
with the only edge of $\partial(F)$ (that is not a half-edge).
As an example,  $M_2$ is shown in
Figure \ref{fig3N} (c).


 \begin{figure}[h]
\setlength{\unitlength}{1.5mm}
     \begin{center}
  \begin{tikzpicture}[scale=0.8]

\coordinate(b1) at  (1, 2);
\coordinate(c1) at  (2, 2);
\coordinate(s1) at  (1.5, 1.5);
\coordinate[label=left:$u$] (e1) at  (0.0, 1);
\coordinate[label=right:$v$](f1) at  (3, 1);

\draw[thick](e1)--(b1)--(c1)--(s1)--(e1)(s1)--(f1)--(e1);
\draw[thick](e1)..controls(1.4, 3.1)..(c1);
\draw[thick](b1)..controls(1.6, 2.8)..(2, 2.5)..controls(2.61, 2.05)..(s1);
 \filldraw [black](b1) circle (3pt)(c1) circle (3pt) (s1) circle (3pt)(e1) circle (3pt)(f1) circle (3pt);
  \node at (1.45, 0.2) {(a)~$M^*$};


\coordinate(a1) at  (5.5, 1);
\coordinate(b1) at  (8, 1);
\coordinate(c1) at  (8, 4);
\coordinate(d1) at  (5.5, 4);

\draw[thick](a1)--(b1)--(c1)--(d1)--(a1);
\draw[thick](a1)..controls(6, 3)..(c1);
\draw[thick](b1)..controls(7.5, 3)..(d1);

 \filldraw [black](a1) circle (3pt)(b1) circle (3pt)(c1) circle (3pt)(d1) circle (3pt);
  \node at (6.6, 0.2) {(b)~$M_1$};
\coordinate(a1) at  (10.5, 1);
\coordinate(b1) at  (13.5, 1);
\coordinate(c1) at  (13.5, 4);
\coordinate(d1) at  (10.5, 4);
\coordinate(f1) at  (11.5, 2);
\coordinate(g1) at  (12.5, 2);
\coordinate(s1) at  (12, 1.5);
\draw[thick](a1)--(b1)--(c1)--(d1)--(a1)
(a1)--(f1)--(g1)--(s1)--(a1)(b1)--(s1);
\draw[thick](a1)..controls(10.9, 3.5)..(c1);
\draw[thick](b1)..controls(13, 3.5)..(d1);
\draw[thick](a1)..controls(11.9, 3.1)..(g1);
\draw[thick](f1)..controls(12.1, 2.8)..(12.5, 2.5)..controls(13.01, 2.05)..(s1);
\filldraw [black](a1) circle (3pt)(b1) circle (3pt)(c1) circle (3pt)(d1) circle (3pt)(f1) circle (3pt)(g1) circle (3pt) (s1) circle (3pt);

  \node at (11.65, 0.2) {(c)~$M_2$};
\
      \end{tikzpicture}
  \caption{{\small~$M^{*}$, $M_{1}$ and $M_{2}$.}}\label{fig3N}
  \end{center} \end{figure}

Note that for each $k\ge 1$,
	$M_{k}$ is a maximal NIC-plane graph without hermits and
	$$
	|V(M_k)|=3k+1, \qquad
		|E(M_{k})|=6+7(k-1)=7k-1.
	$$
It is easy to see that each $M_k$
has exactly $4+2(k-1)=2k+2$ false $3$-faces.
For each $k\ge 1$,
let $M'_k$ be the graph obtained
by inserting a hermit within each
false 3-face of $M_{k}$.

\begin{pro}\label{M'_k-lem}
	For each $k\ge 1$, $M'_k$
is a maximal NIC-plane graph
with
$n$ vertices and $\frac{11}{5}n-\frac{18}{5}$ edges, where $n=5k+3$.
\end{pro}
\proof Since $M_k$ has exactly $2k+2$ false $3$-faces, by the definition
	of $M'_k$,
 $$
 n=|V(M_k)|+2k+2=3k+1+2k+2=5k+3
 $$and
 $$
|E(M'_k)|=|E(M_k)|+2(2k+2)
=7k-1+2(2k+2)
=11k+3=
\frac{11}{5}n-\frac{18}{5}.
$$
\proofend

This shows the tightness of bound of Theorem \ref{theo2N}.

\section{Proofs of the main results} \label{proof}

By Propositions~\ref{H'_k-lem}
and~\ref{M'_k-lem}, in order to prove Theorems~\ref{theo2} and~\ref{theo2N},
it suffices to show that
$|E(G)|\ge \left\lceil \frac 73 n-\frac{14}{3} \right\rceil$
for any maximal IC-plane graph
$G$, and
$|E(G)|\ge \left\lceil \frac {11}{5} n-\frac{18}{5} \right\rceil$
for any maximal NIC-plane graph
$G$. We start with a special case.

\subsection{A special case}
In this subsection, we prove Theorems~\ref{theo2} and~\ref{theo2N} for a special subfamilies of maximal IC-plane resp., maximal NIC-plane graphs.

Let $G$ be either a maximal IC-plane graph or a maximal NIC-plane graph.
Let ${\cal C}^*(G)$ be the set
	of crossings $\alpha$ in $G$
	such that $\al$ is
	incident with four false 4-faces, as shown in Figure~\ref{6}.

 \begin{figure}[h]
\setlength{\unitlength}{1.5mm}
     \begin{center}
  \begin{tikzpicture}[scale=0.8]

\coordinate(a1) at  (0.5, 1);
\coordinate(b1) at  (3.5, 1);
\coordinate(c1) at  (3.5, 4);
\coordinate(d1) at  (0.5, 4);
\coordinate(v1) at  (1.3, 2.5);
\coordinate(v2) at  (2.7, 2.5);
\coordinate(v3) at  (2, 1.8);
\coordinate(v4) at  (2, 3.2);
\draw[thick](a1)--(b1)--(c1)--(d1)--(a1)--(c1)(b1)--(d1)
(a1)--(v1)--(d1)(a1)--(v3)--(b1)(b1)--(v2)--(c1)(c1)--(v4)--(d1);
 \filldraw [black](a1) circle (3pt)(b1) circle (3pt)(c1) circle (3pt)(d1) circle (3pt)(v1) circle (3pt)(v2) circle (3pt)(v3) circle (3pt)(v4) circle (3pt);
 \node at (2, 2.25) {\small{$\alpha$}};\node at (2, 1.3) {\small{$R_1'$}};   \node at (3.2, 2.5) {\small{$R_2'$}};   \node at (2, 3.7) {\small{$R_3'$}};   \node at (0.8, 2.5) {\small{$R_4'$}};   \node at (0.1, 2.5) {\small{$R_5$}};

\coordinate(a2) at  (8.5, 1);
\coordinate(b2) at  (11.5, 1);
\coordinate(c2) at  (11.5, 4);
\coordinate(d2) at  (8.5, 4);
\coordinate(u1) at  (9.3, 2.5);
\coordinate(u2) at  (10.7, 2.5);
\coordinate(u3) at  (10, 1.8);
\coordinate(u4) at  (10, 3.2);
\draw[thick](a2)--(b2)--(c2)--(d2)(a2)--(c2)(b2)--(d2)
(a2)--(u1)--(d2)(a2)--(u3)--(b2)(b2)--(u2)--(c2)(c2)--(u4)--(d2);
 \filldraw [black](a2) circle (3pt)(b2) circle (3pt)(c2) circle (3pt)(d2) circle (3pt)(u1) circle (3pt)(u2) circle (3pt)(u3) circle (3pt)(u4) circle (3pt);
 \draw[thick](d2)..controls(12.5,5)..(12.5, 3)..controls(12.5,0)..(a2);
 \node at (10, 2.25) {\small{$\alpha$}};\node at (10, 1.3) {\small{$R_1'$}};   \node at (11.2, 2.5) {\small{$R_2'$}};   \node at (10, 3.7) {\small{$R_3'$}};   \node at (8.8, 2.5) {\small{$R_4'$}};   \node at (12, 2.5) {\small{$R_5$}};

      \end{tikzpicture}

  (a) \hspace{5 cm} (b)

  \caption{{Crossing $\alpha$ is incident with four false 4-faces.}}\label{6}
  \end{center}
  \end{figure}

	For any $\alpha\in {\cal C}^*(G)$,
	if $R_1, \, R_2, \,  R_3$ and $R_4$ are
	the four $3$-regions in $\R(\al)$,
	then each $R_i$, where $1\le i\le 4$,
	is divided into
	two subregions,
i.e., a false $4$-face
		which are incident with $\alpha$
		and a subregion,
		denoted by $R_i'$,
		which is bounded by three non-crossing edges,
		as shown in Figure~\ref{6}.
	
	Let $\R^*(\al)=\{R'_i: 1\le i\le 4\}$ and let
$r^*(\al)$ denote the number
	of regions in $\R^*(\al)$
	that are not faces of $G$.
If $\C^*(G)\ne \emptyset$, then $G$ contains a subgraph isomorphic to Figure~\ref{6} (a) or (b), implying that
	$|V(G)|\ge 8$.

\begin{pro}\label{pro4-1N}
Let $G$ be either a maximal IC-plane graph or a maximal NIC-plane graph.
	Assume that each face
	in $G$ is of size
	$3$ or $4$ and
	each crossing of $G$ belongs to $\C^*(G)$.
	If $G$ is a maximal IC-plane graph, then
	$|E(G)|\ge \left\lceil \frac 73 n-\frac{14}{3} \right\rceil$. If $G$ is a maximal NIC-plane graph, then
		$|E(G)|\ge \left\lceil \frac {11}{5} n-\frac{18}{5} \right\rceil$.
\end{pro}

\def \G {{\cal G}}

\proof  Let $\G^*$ denote
the set of maximal IC-plane (or maximal NIC-plane) graphs $G$
such that
each face in $G$ is of size
$3$ or $4$ and
each crossing of
$G$ belongs to $\C^*(G)$.

Suppose that Proposition~\ref{pro4-1N}
fails, and $G\in \G^*$ is a counter-example such that its order $n$ is as small as possible.

Clearly, $n\ge 3$.
If $G$ has no crossings, then
$G$ is a maximal plane graph,
implying that $|E(G)|=3n-6
\ge \max\left\{\frac 73 n-\frac{14}{3},
\frac {11}{5} n-\frac{18}{5} \right\}$,
a contradiction to the assumption on $G$.
Thus, $G$ has at least one crossing.
Since each crossing of $G$ belongs to $\C^*(G)$,
$G$ contains a subgraph isomorphic to Figure~\ref{6} (a) or (b),
implying that $n\ge 8$.
If $n=8$, then $G$ is a graph
as shown in Figure~\ref{6}
which has $14$ edges.
However, when $n=8$,
$\frac{7}{3}n-\frac{14}{3}= \frac{11}{5}n-\frac{18}{5}=14$,
a contradiction to the assumption on $G$.
Hence, $n\ge 9$.

Assume that $G$ has $c$ crossings and $h$ hermits. By Property~\ref{propb}~(1) and (2), $G$ has at most $5c$ 4-faces. By Lemma~\ref{lemaN}~(2),
for any 4-face
$F$,
there exist at least three vertices on
$\partial(F)$.
Then, we can obtain a new graph $G^{P}$, which may not be simple,
by adding one edge
to each 4-face in $G$,
so that each face in $G^{P}$ is a 3-face.
Observe that $|E(G^{P})|\leq|E(G)|+5c.$
By Lemma~\ref{lemB}, $|E(G^{P})|=3n-6+c.$
Thus, Claim 1 below holds.

\begin{cla}\label{cl5N}
$|E(G)|\geq  3n-6-4c$.
\end{cla}

\begin{cla}\label{cl6N}
For each $\al\in \C^*(G)$,
$r^*(\al)\le 1$, i.e.,
at most one region in $\R^*(\al)$
is not a face of $G$.
\end{cla}

Suppose that Claim~\ref{cl6N} fails.  Then,
there exists $\al\in \C^*(G)$
with $r^*(\al)\ge 2$.
Assume that
$R'_{1}$ and $R'_2$ are not faces of $G$.

 \begin{figure}[h]
\setlength{\unitlength}{1.5mm}
     \begin{center}
  \begin{tikzpicture}[scale=0.8]

\coordinate(a1) at  (0.5, 1);
\coordinate(b1) at  (3.5, 1);
\coordinate(c1) at  (3.5, 4);
\coordinate(d1) at  (0.5, 4);
\coordinate(v1) at  (1.3, 2.5);
\coordinate(v2) at  (2.7, 2.5);
\coordinate(v3) at  (2, 1.8);
\coordinate(v4) at  (2, 3.2);
\draw[thick](a1)--(b1)--(c1)--(d1)--(a1)--(c1)(b1)--(d1)
(a1)--(v1)--(d1)(a1)--(v3)--(b1)(b1)--(v2)--(c1)(c1)--(v4)--(d1);
 \filldraw [black](a1) circle (3pt)(b1) circle (3pt)(c1) circle (3pt)(d1) circle (3pt)(v1) circle (3pt)(v2) circle (3pt)(v3) circle (3pt)(v4) circle (3pt);
 \node at (2, 1.3) {\small{$R_1'$}};   \node at (3.2, 2.5) {\small{$R_2'$}};   \node at (2, 3.7) {\small{$R_3'$}};   \node at (0.8, 2.5) {\small{$R_4'$}};   \node at (0.1, 2.5) {\small{$R_5$}};
  \node at (2.05, 0.2) {(a)~$G$};

\coordinate(a2) at  (5.5, 1);
\coordinate(b2) at  (8.5, 1);
\coordinate(c2) at  (8.5, 4);
\coordinate(d2) at  (5.5, 4);
\coordinate(s1) at  (6.3, 2.5);
\coordinate(s2) at  (7.7, 2.5);
\coordinate(s3) at  (7, 1.8);
\coordinate(s4) at  (7, 3.2);
\fill[fill=gray, fill opacity=0.2, even odd rule]  (4.8,0.5) rectangle (9.2,4.5)  (a2) rectangle (c2);
\fill[fill=gray, fill opacity=0.2] (a2)--(s1)--(d2)--cycle (d2)--(s4)--(c2)--cycle (c2)--(s2)--(b2)--cycle;
\draw[thick](a2)--(b2)--(c2)--(d2)--(a2)--(c2)(b2)--(d2)
(a2)--(s1)--(d2)(a2)--(s3)--(b2)(b2)--(s2)--(c2)(c2)--(s4)--(d2);
 \filldraw [black](a2) circle (3pt)(b2) circle (3pt)(c2) circle (3pt)(d2) circle (3pt)(s1) circle (3pt)(s2) circle (3pt)(s3) circle (3pt)(s4) circle (3pt);
  \node at (8.2, 2.5) {\small{$R_2'$}};   \node at (7, 3.7) {\small{$R_3'$}};   \node at (5.8, 2.5) {\small{$R_4'$}};   \node at (5.1, 2.5) {\small{$R_5$}};
  \node at (7.05, 0.2){(b)~$G_1$};

\coordinate(a3) at  (10.5, 1);
\coordinate(b3) at  (13.5, 1);
\coordinate(c3) at  (13.5, 4);
\coordinate(d3) at  (10.5, 4);
\coordinate(u1) at  (11.3, 2.5);
\coordinate(u2) at  (12.7, 2.5);
\coordinate(u3) at  (12, 1.8);
\coordinate(u4) at  (12, 3.2);
\fill[fill=gray, fill opacity=0.2] (a3)--(u3)--(b3)--cycle;
\draw[thick](a3)--(b3)--(c3)--(d3)--(a3)--(c3)(b3)--(d3)
(a3)--(u1)--(d3)(a3)--(u3)--(b3)(b3)--(u2)--(c3)(c3)--(u4)--(d3);
 \filldraw [black](a3) circle (3pt)(b3) circle (3pt)(c3) circle (3pt)(d3) circle (3pt)(u1) circle (3pt)(u2) circle (3pt)(u3) circle (3pt)(u4) circle (3pt);
 \node at (12, 1.3) {\small{$R_1'$}};
  \node at (12.05, 0.2) {(c)~$G_2$};
      \end{tikzpicture}
  \caption{{\small $G$, $G_1$ and $G_2$.}}\label{7}
  \end{center}
  \end{figure}

For $1\le i\le 4$,
let $In(R'_i)$ denote the set of
interior vertices of $R'_i$.
Similarly,  let
$In(R_5)$ denote the set of
interior vertices of $R_5$.
We now consider two subgraphs $G_{1}$ and $G_{2}$ of $G$,
where
$G_1=G-In(R'_1)$
(i.e., the subgraph obtained
from $G$ by removing all
vertices in $In(R'_1)$)
and
$G_2=G-\left (In(R_5)\cup
\bigcup\limits_{2\le i\le 4} In(R'_i)\right )$,
as shown in Figure \ref{7} (b) and (c).

We assume that
both $G_{1}$ and $G_{2}$ inherit the original drawing of $G$.
Thus,
both  $G_{1}$ and $G_{2}$
are maximal IC-plane (or maximal NIC-plane),
otherwise, an edge can be added to $G$ without violating the
assumption on $G$.

Observe that both $G_{1}$ and $G_{2}$
belong to $\G^*$ and
 satisfy all conditions of
Proposition~\ref{pro4-1N}.
Since both $R'_1$ and $R'_2$ are
not faces of $G$,
both $G_1$ and $G_2$ have less than $n$ vertices.
By the assumption on $G$,  for $i=1,2$,
$|E(G_{i})|\geq c|V(G_i)|-d,$
where $c=\frac{7}{3},d=\frac{14}{3}$,
if $G_i$ is a maximal IC-plane graph,
or $c=\frac{11}{5},d=\frac{18}{5}$
otherwise.
It can be verified that $d=8c-14$.

Note that  $|V(G_{1})|+|V(G_{2})|=n+8$ and $|E(G_{1})|+|E(G_{2})|=|E(G)|+14$.
Thus,
\begin{eqnarray*}
	|E(G)|&=&|E(G_{1})|+|E(G_{2})|-14\\
	&\geq& \Big(c|V(G_1)|-d\Big)+\Big(c|V(G_2)|-d\Big)-14\\
	&=&cn+8c-2d-14=cn-d,
\end{eqnarray*}
a contradiction to the assumption on $G$. Hence, Claim~\ref{cl6N} holds.

By Claim~\ref{cl6N}, for each crossing
$\al\in \C^*(G)$,  $r^*(\al)\le 1$,
implying that
each crossing $\al$ in $G$
corresponds to at least three different hermits.
Thus, $h\geq 3c$.

If $G$ is a maximal IC-plane graph, by Lemma~\ref{lemc},
$c\leq \frac{n-3c}{4}$,
implying that
$c\leq \frac{n}{7}$.
By Claim~\ref{cl5N},
$|E(G)|\geq  3n-6-4c$.
Thus,
$|E(G)|\geq \frac{17}{7}n-6$,
implying that
$$
|E(G)|\geq
	\left	\lceil \frac{17}{7}n-6
	\right \rceil
	\geq  \frac{7}{3}n-\frac{14}{3}
$$
for all $n\geq 10$.
Since $G$ is a counter-example,
we have $n=9$.
Then, $In(R'_1)=\{z\}$
for some vertex $z$.
By maximality of $G$,
we have $d_G(z)=3$,
implying that
$|E(G)|=14+3\geq  \frac{7}{3}n-\frac{14}{3}$
for $n=9$,
a contradiction.

Now, we consider the case
that $G$ is a
maximal NIC-plane graph.
By Lemma~\ref{lemcN}, $c\leq \frac{|E(G)|-6c}{6}$,
implying that
$c\leq \frac{|E(G)|}{12}$.
By Claim~\ref{cl5N},
$|E(G)|\geq  3n-6-4c$.
Thus,
$|E(G)|\geq \frac{9}{4}n-\frac{9}{2}$,
implying that
$$
|E(G)|\geq
	\left	\lceil \frac{9}{4}n-\frac{9}{2}
	\right \rceil
	\geq \left	\lceil \frac{11}{5}n-\frac{18}{5}\right \rceil
$$
for all $n\geq 11$ with $n\ne 14$.

Since $G$ is a counter-example,
we have $n\in \{9,10,14\}$ and
$|E(G)|<\frac{11}{5}n-\frac{18}{5}$.
As $|E(G)|\geq  3n-6-4c$,
we have
$
3n-6-4c<\frac{11}{5}n-\frac{18}{5},
$
implying that $c>\frac {n-3}5$.
However,  for each $n\in \{9,10,14\}$,
as $|E(G)|\ge 12 c $, we have
$$
	|E(G)|\ge 12 c \ge 12
	\left \lceil \frac {n-3}5 \right \rceil
	> \left \lceil\frac{11}{5}n-\frac{18}{5} \right \rceil,
$$
a contradiction.
Thus,  $G$ is not a counter-example
to the conclusion of Proposition~\ref{pro4-1N}.
Proposition~\ref{pro4-1N} holds.
\proofend

\vspace{0.2 cm}

\subsection{General case}

We are ready to prove Theorems~\ref{theo2} and~\ref{theo2N}.

\noindent \textbf{\emph{Proof of Theorem~\ref{theo2}.}}
By Proposition~\ref{H'_k-lem},
it suffices to show that
$|E(G)|\ge
\left\lceil \frac{7}{3}n-\frac{14}{3} \right \rceil$
for any maximal IC-plane graph $G$.
Suppose that
$G$ is a maximal IC-plane graph
of order $n$ such that $n$ is as small as possible and
$|E(G)|<\frac{7}{3}n
-\frac{14}{3}$.
Clearly, $n\ge 5$ and we can assume that
$G$ has at least one crossing.

Assume that $G$ has exactly
$t$  false 3-faces.
We construct a new graph $\widetilde{G}$ by adding
a new vertex $u$ inside each
false 3-face of $G$
and adding a new edge
joining $u$ to both vertices
on $\partial(F)$.
Then, we have $|V(\widetilde{G})|=n+t$ and $|E(\widetilde{G})|=|E(G)|+2t$.

Observe that $\widetilde{G}$ is a maximal IC-plane
graph and satisfies all conditions of
Proposition~\ref{pro4-1N}.
By Proposition~\ref{pro4-1N},
$$
	|E(\widetilde{G})|
	\ge  \frac{7}{3}|V(\widetilde{G})|-\frac{14}{3}=\frac{7}{3}(n+t)-\frac{14}{3}.
$$
It follows that $$
|E(G)|=|E(\widetilde{G})|-2t
\ge \frac{7}{3}n+\frac t3-\frac{14}{3}
\geq \frac{7}{3}n-\frac{14}{3},
$$
a contradiction to the assumption on $G$. Hence,
$|E(G)|\ge \left \lceil \frac{7}{3}n-\frac{14}{3} \right \rceil$
for each maximal IC-plane graph $G$.
\proofend

\noindent \textbf{\emph{Proof of Theorem~\ref{theo2N}.}}

By Proposition~\ref{M'_k-lem},
	it suffices to show that
	$|E(G)|\ge
	\left\lceil \frac{11}{5}n-\frac{18}{5}\right\rceil$
	for any maximal NIC-plane graph $G$.
Suppose that
$G$ is a maximal NIC-plane graph
of order $n$ such that $n$ is as small as possible and
$|E(G)|<\frac{11}{5}n-\frac{18}{5}$. Clearly, $n\ge 5$, and we can assume that $G$ has at least one crossing.

\begin{cla}\label{cl8N_1}
There are no maximal NIC-plane
graphs $G_1$ and $G_2$
such that
$|V(G_i)|<|V(G)|$
for both $i=1, \, 2$,
$|V(G_1)|+|V(G_2)|=n+3$
and $|E(G_1)|+|E(G_2)|=|E(G)|+3$.
\end{cla}

Suppose that such maximal NIC-plane graphs $G_1$ and
$G_2$ exist.
Then, by the assumption on $G$,
$|E(G_{i})| \geq \frac{11}{5}|V(G_i)|-\frac{18}{5}$
for both $i=1,2$.
Since  $|V(G_{1})|+|V(G_{2})|=n+3$
and $|E(G_1)|+|E(G_2)|=|E(G)|+3$,
\eqn{th2N-eq2}
{
	|E(G)|&=&
	|E(G_{1})|+|E(G_{2})|-3
	\nonumber \\
	&\geq& \Big(\frac{11}{5}|V(G_1)|
	-\frac{18}{5}\Big)
	+\Big(\frac{11}{5}|V(G_2)|
	-\frac{18}{5}\Big)-3
	\nonumber \\
	&=&\frac{11}{5}n-\frac{18}{5},
}
a contradiction to the assumption on $G$.
Thus, Claim~\ref{cl8N_1} holds.

\begin{cla}\label{cl8N}
$G$ has no $5$-faces.
\end{cla}

Suppose that Claim~\ref{cl8N} fails, and $F$ is a
$5$-face in $G$.
By Lemma~\ref{lemaN} (3),
$\alpha(F)$ consists of one edge $ab$
and four half-edges $\alpha_1b, \, \alpha_1c, \,  \alpha_2a$
and $\alpha_2c$,
where $\alpha_1$ and $\alpha_2$ are crossings,
implying that
$G$ is as shown in Figure \ref{fig5N} (a).
\begin{figure} [h!]
  \centering
\begin{tikzpicture}[scale=0.39]
	\begin{pgfonlayer}{nodelayer}

		\node [style=blacknode] (7) at (-9.325, 2.35) {};
		\node [style=blacknode] (8) at (-8.25, 1.275) {};
		\node [style=blacknode] (10) at (-3.5, 0) {};
		\node [style=blacknode] (12) at (3.5, 0) {};
		\node [style=blacknode] (13) at (0, 3.5) {};

		\node [style=blacknode] (16) at (1.175, 2.35) {};
		\node [style=blacknode] (17) at (2.25, 1.275) {};
		\node [style=blacknode] (18) at (7, 7) {};
		\node [style=blacknode] (19) at (7, 0) {};
		\node [style=blacknode] (20) at (14, 7) {};
		\node [style=blacknode] (21) at (14, 0) {};
		\node [style=blacknode] (22) at (10.5, 3.5) {};
		\node [style=none] (23) at (-15, 8) {};
		\node [style=none] (24) at (-6, 8) {};
		\node [style=none] (25) at (-6, -1) {};
		\node [style=none] (26) at (-15, -1) {};
		\node [style=none] (27) at (6, 8) {};
		\node [style=none] (28) at (15, 8) {};
		\node [style=none] (29) at (15, -1) {};
		\node [style=none] (30) at (6, -1) {};
		\node [style=none] (31) at (10.5, 4.25) {$b$};

		\node [style=none] (32) at (10.5, 5.311) {$\alpha_2$};

		\node [style=none] (34) at (2.65, 2.10) {$\alpha_1$};
		\node [style=none] (35) at (0, 4.2) {$b$};

		\node [style=none] (36) at (-4, 0.25) {$a$};
		\node [style=none] (37) at (4, 0.25) {$c$};

		\node [style=none] (38) at (-10.5, 1) {$R$};
		\node [style=none] (39) at (-8, 2.30) {$\alpha_1$};
		\node [style=none] (40) at (-10.5, 4.2) {$b$};
		\node [style=none] (42) at (-14.5, 0.25) {$a$};
		\node [style=none] (43) at (-6.5, 0.25) {$c$};

		\node [style=blacknode] (49) at (-14, 7) {};
		\node [style=blacknode] (50) at (-14, 0) {};
		\node [style=blacknode] (51) at (-7, 7) {};
		\node [style=blacknode] (52) at (-7, 0) {};
		\node [style=blacknode] (53) at (-10.5, 3.5) {};
				\node [style=none] (55) at (-10.5, 5.311) {$\alpha_2$};

\node [style=none] (59) at (-10.5, -2) {(a) $G$};
\node [style=none] (60) at (0, -2) {(b) $G_{1}$};
\node [style=none] (60) at (10.5, -2) {(c) $G_{2}$};
		\node [style=none] (56) at (6.5, 0.25) {$a$};
		\node [style=none] (57) at (14.5, 0.25) {$c$};
		\node [style=none] (58) at (-12, 4) {\textbf{$F$}};
		\node [style=none] (60) at (0, 1.25) {\textbf{$R$}};
	\end{pgfonlayer}
	\begin{pgfonlayer}{edgelayer}

\fill[fill=gray!20] (23.center) to (49.center) to (51.center) to (24.center);
\fill[fill=gray!20] (25.center) to (52.center) to (51.center) to (24.center);
\fill[fill=gray!20] (23.center) to (49.center) to (51.center) to (24.center);
\fill[fill=gray!20] (23.center) to (49.center) to (50.center) to (26.center);
\fill[fill=gray!20] (25.center) to (52.center) to (50.center) to (26.center);
\fill[fill=gray!20] (49.center) to (50.center) [in=-165, out=83, looseness=1.25] to (51.center);
\fill[fill=gray!20] (51.center) to (49.center) [in=97, out=-15, looseness=1.25] to (52.center);

\fill[fill=gray!20] (53.center) [bend left] to (50.center) [bend right] to  (52.center);
\fill[fill=gray!20] (53.center) [bend left=45] to (8.center)to (52.center);
\fill[fill=gray!20] (7.center) [bend left] to (52.center) to (53.center);

\fill[fill=gray!20] (12.center) to (13.center) [bend left] to (10.center);

\fill[fill=gray!20] (13.center) [bend left=45] to (17.center) to (12.center);
\fill[fill=gray!20] (16.center) [bend left]  to (12.center) to (13.center);

\fill[fill=gray!20] (27.center) to (18.center) to (20.center) to (28.center);
\fill[fill=gray!20] (29.center) to (21.center) to (20.center) to (28.center);
\fill[fill=gray!20] (29.center) to (21.center) to (19.center) to (30.center);
\fill[fill=gray!20] (30.center) to (19.center) to (18.center) to (27.center);

\fill[fill=gray!20] (18.center) to (19.center) [in=-165, out=83, looseness=1.25] to (20.center);
\fill[fill=gray!20] (20.center) to (18.center) [in=97, out=-15, looseness=1.25] to (21.center);

		\draw [style=blackedge_thick, in=180, out=0] (49.center) to (51.center);
		\draw [style=blackedge_thick, in=-165, out=83, looseness=1.25] (50.center) to (51.center);
		\draw [style=blackedge_thick, in=97, out=-15, looseness=1.25] (49.center) to (52.center);

		\draw [style=blackedge_thick, bend left] (53.center) to (50.center);
		\draw [style=blackedge_thick, bend left=45] (53.center) to (8.center);
		\draw [style=blackedge_thick, bend left] (7.center) to (52.center);
		\draw [style=blackedge_thick, bend right] (53.center) to (52.center);
		\draw [style=blackedge_thick] (53.center) to (52.center);
		\draw [style=blackedge_thick] (51.center) to (52.center);
		\draw [style=blackedge_thick, in=180, out=0] (50.center) to (52.center);
		\draw [style=blackedge_thick, in=90, out=-90] (49.center) to (50.center);
		\draw [style=blackedge_thick, in=180, out=0] (10.center) to (12.center);
		\draw [style=blackedge_thick, bend left] (13.center) to (10.center);
		\draw [style=blackedge_thick, bend right] (13.center) to (12.center);

		\draw [style=blackedge_thick] (13.center) to (12.center);
		\draw [style=blackedge_thick, bend left=45] (13.center) to (17.center);
		\draw [style=blackedge_thick, bend left] (16.center) to (12.center);
		\draw [style=blackedge_thick, in=180, out=0] (18.center) to (20.center);
		\draw [style=blackedge_thick] (20.center) to (21.center);
		\draw [style=blackedge_thick, in=180, out=0] (19.center) to (21.center);
		\draw [style=blackedge_thick] (18.center) to (19.center);
		\draw [style=blackedge_thick, in=-165, out=83, looseness=1.25] (19.center) to (20.center);
		\draw [style=blackedge_thick, in=97, out=-15, looseness=1.25] (18.center) to (21.center);
		\draw [style=blackedge_thick, bend left] (22.center) to (19.center);
		\draw [style=blackedge_thick, bend right] (22.center) to (21.center);
		\draw [style=blackedge_thick, in=180, out=0] (49.center) to (51.center);
		\draw [style=blackedge_thick] (51.center) to (52.center);
		\draw [style=blackedge_thick, in=180, out=0] (50.center) to (52.center);
		\draw [style=blackedge_thick] (49.center) to (50.center);
		\draw [style=blackedge_thick, bend right] (53.center) to (52.center);
	\end{pgfonlayer}
\end{tikzpicture}

  \caption{$G$, $G_{1}$ and  $G_{2}$.}\label{fig5N}
\end{figure}
Let $G_1$ be the subgraph of $G$ induced by $V_1$,
where $V_1$ is
the set of vertices which are within the finite region bounded by $ab\alpha_1ca$,
including vertices on the boundary,
and let  $G_2$ be the subgraph of $G$ induced by $V_2$,
where $V_2=(V(G)\setminus V_1)\cup \{a, \, b, \, c\}$.
We may assume that
both $G_1$ and $G_2$ inherit the drawing of $G$,
as shown in Figure \ref{fig5N} (b) and (c).

Clearly, $G_{2}$ is maximal NIC-plane.
We claim that
$G_{1}$ is also maximal NIC-plane. Indeed,
otherwise,
there is either an edge
passing through the edge $ac$
and
connecting vertex $b$ to some
vertex within region $R$,
or an edge passing through the edge $ab$ and
connecting vertex $c$ to some vertex
within region $R$.
However, both
situations violate the
assumption on $G$,
a contradiction.

Note that both $G_{1}$ and $G_{2}$
have fewer vertices than $G$,
$|V(G_1)|+|V(G_2)|=|V(G)|+3=n+3$
and
$|E(G_1)|+|E(G_2)|=|E(G)|+3$.
However,
by Claim~\ref{cl8N_1},
such maximal NIC-plane graphs $G_1$ and $G_2$
do not exist.
Hence, Claim~\ref{cl8N}
holds.

\begin{cla}\label{cl9N}
$G$ has no $6$-faces.
\end{cla}

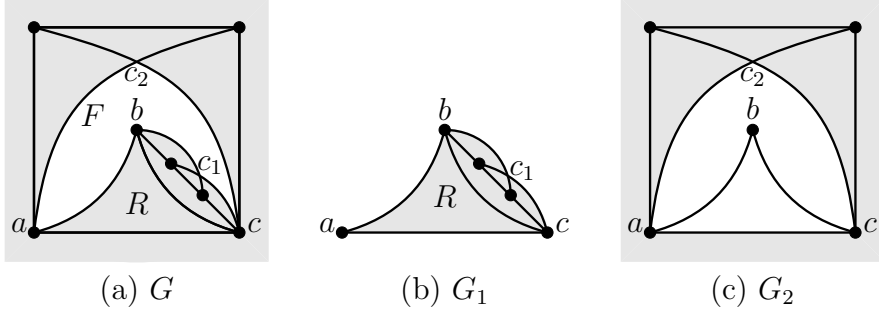
\begin{figure} [h!]
  \centering
\begin{tikzpicture}[scale=0.4]
	\begin{pgfonlayer}{nodelayer}
		\node [style=blacknode] (6) at (-12.75, 1.25) {};
		\node [style=blacknode] (7) at (-9.325, 2.35) {};
		\node [style=blacknode] (8) at (-8.25, 1.275) {};
		\node [style=blacknode] (10) at (-3.5, 0) {};
		\node [style=blacknode] (12) at (3.5, 0) {};
		\node [style=blacknode] (13) at (0, 3.5) {};
		\node [style=blacknode] (14) at (-1.25, 2.35) {};
		\node [style=blacknode] (15) at (-2.25, 1.25) {};
		\node [style=blacknode] (16) at (1.175, 2.35) {};
		\node [style=blacknode] (17) at (2.25, 1.275) {};
		\node [style=blacknode] (18) at (7, 7) {};
		\node [style=blacknode] (19) at (7, 0) {};
		\node [style=blacknode] (20) at (14, 7) {};
		\node [style=blacknode] (21) at (14, 0) {};
		\node [style=blacknode] (22) at (10.5, 3.5) {};
		\node [style=none] (23) at (-15, 8) {};
		\node [style=none] (24) at (-6, 8) {};
		\node [style=none] (25) at (-6, -1) {};
		\node [style=none] (26) at (-15, -1) {};
		\node [style=none] (27) at (6, 8) {};
		\node [style=none] (28) at (15, 8) {};
		\node [style=none] (29) at (15, -1) {};
		\node [style=none] (30) at (6, -1) {};
		\node [style=none] (31) at (10.5, 4) {$c$};

		\node [style=none] (32) at (10.5, 5.311) {$\alpha_3$};
		\node [style=none] (33) at (-2.5, 2.45) {$\alpha_1$};
		\node [style=none] (34) at (2.65, 2.15) {$\alpha_2$};
		\node [style=none] (35) at (0, 4) {$c$};

		\node [style=none] (36) at (-4, 0.25) {$a$};
		\node [style=none] (37) at (4, 0.25) {$b$};

		\node [style=none] (38) at (-12.9, 2.45) {$\alpha_1$};
		\node [style=none] (39) at (-8, 2.35) {$\alpha_2$};
		\node [style=none] (40) at (-10.5, 4) {$c$};
		\node [style=none] (42) at (-14.5, 0.25) {$a$};
		\node [style=none] (43) at (-6.5, 0.25) {$b$};
		\node [style=blacknode] (48) at (-11.75, 2.35) {};
		\node [style=blacknode] (49) at (-14, 7) {};
		\node [style=blacknode] (50) at (-14, 0) {};
		\node [style=blacknode] (51) at (-7, 7) {};
		\node [style=blacknode] (52) at (-7, 0) {};
		\node [style=blacknode] (53) at (-10.5, 3.5) {};
				\node [style=none] (55) at (-10.5, 5.311) {$\alpha_3$};

\node [style=none] (59) at (-10.5, -2) {(a) $G$};
\node [style=none] (60) at (0, -2) {(b) $G_{1}$};
\node [style=none] (60) at (10.5, -2) {(c) $G_{2}$};
		\node [style=none] (56) at (6.5, 0.25) {$a$};
		\node [style=none] (57) at (14.5, 0.25) {$b$};
		\node [style=none] (58) at (-11.55, 4.35) {\textbf{$F$}};
	\end{pgfonlayer}
	\begin{pgfonlayer}{edgelayer}

\fill[fill=gray!20] (23.center) to (49.center) to (51.center) to (24.center);
\fill[fill=gray!20] (25.center) to (52.center) to (51.center) to (24.center);
\fill[fill=gray!20] (23.center) to (49.center) to (51.center) to (24.center);
\fill[fill=gray!20] (23.center) to (49.center) to (50.center) to (26.center);
\fill[fill=gray!20] (25.center) to (52.center) to (50.center) to (26.center);
\fill[fill=gray!20] (49.center) to (50.center) [in=-165, out=83, looseness=1.25] to (51.center);
\fill[fill=gray!20] (51.center) to (49.center) [in=97, out=-15, looseness=1.25] to (52.center);
\fill[fill=gray!20] (53.center) [bend right=45] to (6.center) to (50.center);
\fill[fill=gray!20] (48.center) [bend right, looseness=1.25] to (50.center) to (53.center);
\fill[fill=gray!20] (53.center)  to (50.center) to (52.center);
\fill[fill=gray!20] (53.center) [bend left=45] to (8.center)to (52.center);
\fill[fill=gray!20] (7.center) [bend left] to (52.center) to (53.center);

\fill[fill=gray!20] (13.center) to (10.center) to (12.center);
\fill[fill=gray!20] (13.center) [bend right=45] to (15.center) to (10.center);
\fill[fill=gray!20] (14.center) [bend right, looseness=1.25] to (10.center) to (13.center);
\fill[fill=gray!20] (13.center) [bend left=45] to (17.center) to (12.center);
\fill[fill=gray!20] (16.center) [bend left]  to (12.center) to (13.center);

\fill[fill=gray!20] (27.center) to (18.center) to (20.center) to (28.center);
\fill[fill=gray!20] (29.center) to (21.center) to (20.center) to (28.center);
\fill[fill=gray!20] (29.center) to (21.center) to (19.center) to (30.center);
\fill[fill=gray!20] (30.center) to (19.center) to (18.center) to (27.center);

\fill[fill=gray!20] (18.center) to (19.center) [in=-165, out=83, looseness=1.25] to (20.center);
\fill[fill=gray!20] (20.center) to (18.center) [in=97, out=-15, looseness=1.25] to (21.center);

		\draw [style=blackedge_thick, in=180, out=0] (49.center) to (51.center);
		\draw [style=blackedge_thick, in=-165, out=83, looseness=1.25] (50.center) to (51.center);
		\draw [style=blackedge_thick, in=97, out=-15, looseness=1.25] (49.center) to (52.center);
		\draw [style=blackedge_thick, bend right=45] (53.center) to (6.center);
		\draw [style=blackedge_thick, bend right, looseness=1.25] (48.center) to (50.center);
		\draw [style=blackedge_thick] (53.center) to (50.center);
		\draw [style=blackedge_thick, bend left] (53.center) to (50.center);
		\draw [style=blackedge_thick, bend left=45] (53.center) to (8.center);
		\draw [style=blackedge_thick, bend left] (7.center) to (52.center);
		\draw [style=blackedge_thick, bend right] (53.center) to (52.center);
		\draw [style=blackedge_thick] (53.center) to (52.center);
		\draw [style=blackedge_thick] (51.center) to (52.center);
		\draw [style=blackedge_thick, in=180, out=0] (50.center) to (52.center);
		\draw [style=blackedge_thick, in=90, out=-90] (49.center) to (50.center);
		\draw [style=blackedge_thick, in=180, out=0] (10.center) to (12.center);
		\draw [style=blackedge_thick, bend left] (13.center) to (10.center);
		\draw [style=blackedge_thick, bend right] (13.center) to (12.center);
		\draw [style=blackedge_thick] (13.center) to (10.center);
		\draw [style=blackedge_thick] (13.center) to (12.center);
		\draw [style=blackedge_thick, bend right=45] (13.center) to (15.center);
		\draw [style=blackedge_thick, bend right, looseness=1.25] (14.center) to (10.center);
		\draw [style=blackedge_thick, bend left=45] (13.center) to (17.center);
		\draw [style=blackedge_thick, bend left] (16.center) to (12.center);
		\draw [style=blackedge_thick, in=180, out=0] (18.center) to (20.center);
		\draw [style=blackedge_thick] (20.center) to (21.center);
		\draw [style=blackedge_thick, in=180, out=0] (19.center) to (21.center);
		\draw [style=blackedge_thick] (18.center) to (19.center);
		\draw [style=blackedge_thick, in=-165, out=83, looseness=1.25] (19.center) to (20.center);
		\draw [style=blackedge_thick, in=97, out=-15, looseness=1.25] (18.center) to (21.center);
		\draw [style=blackedge_thick, bend left] (22.center) to (19.center);
		\draw [style=blackedge_thick, bend right] (22.center) to (21.center);
		\draw [style=blackedge_thick, in=180, out=0] (49.center) to (51.center);
		\draw [style=blackedge_thick] (51.center) to (52.center);
		\draw [style=blackedge_thick, in=180, out=0] (50.center) to (52.center);
		\draw [style=blackedge_thick] (49.center) to (50.center);
		\draw [style=blackedge_thick, bend right] (53.center) to (52.center);
	\end{pgfonlayer}
\end{tikzpicture}

  \caption{$G$, $G_{1}$ and  $G_{2}$.}\label{fig6N}
\end{figure}

Suppose that Claim~\ref{cl9N} fails, and $F$
is a $6$-face in $G$.
By Lemma~\ref{lemaN} (3),
$\partial(F)$ consists of six half edges
$\alpha_1a, \; \alpha_1c, \, \alpha_2c, \, \alpha_2b, \, \alpha_3a$ and $\alpha_3b$,
where $\alpha_1,\alpha_2$ and $\alpha_3$ are
crossings in $G$.
Thus, $G$ is as shown in
Figure \ref{fig6N} (a).

Let $G_1$ be the subgraph of $G$ induced by $V_1$,
where $V_1$ is
the set of vertices which are within the finite region bounded by $a\alpha_1c\alpha_2ba$,
including the vertices on the boundary,
and let  $G_2$ be the subgraph of $G$ induced by $V_2$,
where $V_2=(V(G)\setminus V_1)\cup \{a, \, b, \, c\}$.
We may assume that
both $G_1$ and $G_2$ inherit the drawing of $G$,
as shown in Figure \ref{fig6N} (b) and (c).

It can be verified that both
$G_{1}$ and $G_{2}$ are maximal NIC-plane graphs which have  fewer vertices than $G$.
Also note that  $|V(G_{1})|+|V(G_{2})|=n+3$ and $|E(G_{1})|+|E(G_{2})|=|E(G)|
+3$.
However, by Claim~\ref{cl8N_1},
such maximal NIC-plane graphs $G_1$ and $G_2$ do not exist.
Hence, Claim~\ref{cl9N} holds.

By Claims~\ref{cl8N} and \ref{cl9N},
each face in $G$ is of size $3$ or $4$.
Assume that  $G$ has exactly
$t$  false 3-faces.
We construct a new graph $G^{\prime}$ by adding
a new vertex $u$ inside each
false 3-face $F$ of $G$
and adding a new edge
joining $u$ to  each vertex on $\partial(F)$.
Then  $|V(G^{\prime})|=|V(G)|+t=n+t$ and $|E(G^{\prime})|=|E(G)|+2t$.

Note that
$G^{\prime}$ is a maximal NIC-plane graph
and satisfies all conditions of
Proposition~\ref{pro4-1N}.
By Proposition~\ref{pro4-1N},
$$
	|E(G^{\prime})|
	\ge  \frac{11}{5}|V(G^{\prime})|-\frac{18}{5}=\frac{11}{5}(n+t)-\frac{18}{5}.
$$
It follows that $$
|E(G)|=|E(G^{\prime})|-2t
\ge \frac{11}{5}n+\frac t5-\frac{18}{5}
\geq \frac{11}{5}n-\frac{18}{5},
$$
a contradiction to the assumption on $G$. Hence,
$|E(G)|\ge \frac{11}{5}n-\frac{18}{5}$
for each maximal NIC-plane graph $G$.
\proofend

\textbf{Acknowledgments}

The authors are very grateful to the anonymous referees for many comments and suggestions.

\textbf{Declarations}

We claim that there is no conflict of interest in our paper. No data was used for the research described in the article.

\end{document}